\documentclass[10pt,letterpaper]{article}
\makeatletter
\usepackage{amscd}
\usepackage{amsmath}
\numberwithin{equation}{section}
\usepackage{latexsym}
\usepackage{amsfonts}
\usepackage{amssymb}
\usepackage{amsthm}
\usepackage{graphicx}
\usepackage[dvipdfm,
            pdfstartview=FitH,
            CJKbookmarks=true,
            bookmarksnumbered=true,
            bookmarksopen=true,
            colorlinks=true, 
            pdfborder=001,   
            citecolor=magenta,
            linkcolor=blue,
            ]{hyperref}       

\newtheorem{thm}{\textbf Theorem}[section]
\newtheorem{lem}{\textbf Lemma}[section]
\newtheorem{rem}{\textbf Remark}[section]

\newcommand{\md}{\mbox{d}}
\newcommand{\be}{\begin{eqnarray}}
\newcommand{\ee}{\end{eqnarray}}
\newcommand{\mr}{\mathbb{R}}

\newcommand{\ep}{\varepsilon}

\newcommand{\bes}{\begin{eqnarray*}}
\newcommand{\ees}{\end{eqnarray*}}

\begin{document}
\begin{titlepage}
\title{\bf Global large solution for the tropical climate model with diffusion}
\author{ Xia Chen, Baoquan Yuan\thanks{Corresponding Author: B. Yuan}\ , Ying Zhang
       \\ School of Mathematics and Information Science,
       \\ Henan Polytechnic University,  Henan,  454000,  China.\\
        (cx18332162035@163.com, bqyuan@hpu.edu.cn, yz503576344@163.com)
        \date{}}
\end{titlepage}
\maketitle
\begin{abstract}
This paper studies the d-dimensional (d=2,3) tropical climate model with only the dissipation of the first baroclinic model of the velocity ($-\eta\Delta v$). By choosing a class of special initial data $(u_0,v_0,\theta_0)$ whose $H^s(\mr^d)$ norm can be arbitrarily large, we obtain the global smooth solution of d-dimensional tropical climate model.
\end{abstract}

\vspace{.2in} {\bf Key words:}\quad Tropical climate model, global smooth solution, large initial data.

\vspace{.2in} {\bf MSC(2000):} 35Q35, 76D03, 86A10.



\section{Introduction}
\setcounter{equation}{0}
\vskip .1in

This paper focuses on the d-dimensional (d=2,3) tropical climate model (TCM) given by
\be\label{1.1}
\begin{cases}
{ \begin{array}{ll} \partial_{t}u+(u\cdot\nabla) u+\nu u+\nabla p+\textup{div}(v\otimes v)=0, \ \ \ x\in\mr^d, t>0,\\
\partial_{t}v+(u\cdot\nabla) v+(v\cdot\nabla) u-\eta\Delta v+\nabla\theta=0,\\
\partial _{t}\theta+(u\cdot\nabla)\theta+\lambda\theta+\nabla\cdot v=0,\\
\nabla\cdot u=0,\\
u(x,0)=u_{0}(x), v(x,0)=v_{0}(x), \theta(x,0)=\theta_{0}(x),
 \end{array} }
\end{cases}
\ee
where vector fields $u(x,t)$ and $v(x,t)$ represent the barotropic model and the first baroclinic model of the velocity field, respectively. $p=p(x,t)$ and $\theta=\theta(x,t)$ stand for the scalar pressure and scalar temperature. The parameters $\nu$, $\eta$, $\lambda$ are non-negative.

When $\nu=\eta=\lambda=0$, the system (\ref{1.1}) reduces to the original tropical climate model derived by Frierson, Majda and Pauluis \cite{FMP2004} by performing a Galerkin truncation to the hydrostatic Boussinesq. The issues on well-posedness, global regularity and blow-up criterion of TCM have attracted considerable attention. Let us briefly recall some progress. By introducing a new velocity variable, Li-Titi \cite{LT2016} constructed a unique global strong solution for the system (\ref{1.1}) with $-\Delta u$, $\nu=\lambda=0$, $\eta=1$ and $d=2$. Focused on the same version as the system in \cite{LT2016}, Li-Zhai-Yin \cite{LZY2019} established the global well-posedness in negative-order Besov spaces with small initial data. Ye \cite{Y2020} obtained the global regularity result of the 2-dimensional zero thermal diffusion TCM with fractional dissipation, given by $(-\Delta)^\alpha u$ and $(-\Delta)^\beta v$ where $\alpha+\beta\geq 2$ ($1<\alpha<2$). To further have a complete view of current studies on TCM, readers can see \cite{DWY2019, DWWZ2019, MJW2016, W2016, W2019, WZP2020, YT2020, Z2018, Z2020}.

When $\theta=0$, the system (\ref{1.1}) reduces to the MHD-like equations. In \cite{LYY2019}, Li-Yang-Yu established a unique global large solution of the 2-dimensional MHD equations with only damping of velocity and magnetic. Later, Dai-Tan-Wu in \cite{DTW2019} introduced an effective approach to construct a class of global large solutions to the $d$-dimensional $(d=2,3)$ MHD equations with any fractional dissipation. Naturally, we consider that whether or not the system (\ref{1.1}) with a class of special large initial data can generate a unique global solution. Motivated by \cite{DTW2019}, we established a class of global solutions of system (\ref{1.1}) by choosing a class of special initial data $(u_0,v_0,\theta_0)$ whose $H^s (s>1+\frac{d}{2})$ norm can be arbitrarily large.

\vspace{.2in}
There are some differences between the 2-dimensional and 3-dimensional cases, so we split our results into two cases. Before presenting our main result of 3-dimensional case, we define a vector field $m_{0}\in C^{\infty}_{0}(\mr^3)$ and two scalar functions $n_{0}, r_0\in C^{\infty}_{0}(\mr^3)$ which are used to construct the large initial data of system (\ref{1.1})
\begin{align}
&\widehat{m_{0}}(\xi)=\Big(\varepsilon^{-1}
\log\log\frac{1}{\varepsilon}\Big)\chi(\xi),&&\xi\in\mr^3,\label{1.2}\\
&\widehat{n_{0}}(\xi)=\widehat{r_{0}}(\xi)=\Big(\varepsilon^{-1}
\Big(\log\log\frac{1}{\varepsilon}\Big)^{\frac{1}{2}}\Big)\overline{\chi}(\xi),&&\xi\in\mr^3,\label{1.3}
\end{align}
where $0<\varepsilon<1$ is a small parameter depending on $\nu$, $\eta$ and $\lambda$. The smooth cutoff functions $\chi=(\chi_1,\chi_2,\chi_3)$ and $\overline{\chi}$ satisfy
\begin{center}
supp $\chi_i,\overline{\chi}\subset\mathcal{C}$,$ \ $ $\chi_i,\overline{\chi}\in[0,1]$ $ \ $ and $ \ $ $\chi_i,\overline{\chi}=1$ $ \ $ on $ \ $ $\mathcal{C}_{1}$ \ \ $(i=1,2,3)$,
\end{center}
where $\mathcal{C}$ and $\mathcal{C}_1$ denote the annuli,
\begin{align*}
&\mathcal{C}:=\Big\{\xi\in\mr^3\Big{|}\ |\xi_{i}+\xi_{j}|\leq\varepsilon,i,j=1,2,3,i\neq j,1\leq|\xi|\leq 2\Big\},\\
&\mathcal{C}_{1}:=\Big\{\xi\in\mr^3\Big{|}\ |\xi_{i}+\xi_{j}|\leq\varepsilon,i,j=1,2,3,i\neq j,\frac{4}{3}\leq|\xi|\leq\frac{5}{3}\Big\}.
\end{align*}


Now we state our result for 3-dimensional case as follows.

\begin{thm}\label{thm1.1}
Assume $m_{0},n_{0},r_0$ are defined by (\ref{1.2}) and (\ref{1.3}). Consider the initial data in the TCM system (\ref{1.1}) which fulfills $\textup{div}w_{0}=0$ and
\begin{align*}
u_{0}=U_{0}+w_{0},\ \ \ v_{0}=V_{0}+z_{0}\ \ \ and \ \ \ \theta_0=\Theta_0+\psi_0,
\end{align*}
where
\begin{align}\label{1.4}
U_{0}=\nabla\times m_{0},\ \ \ V_{0}=(n_{0},n_{0},n_{0})\ \ \ and\ \ \ \Theta_0=r_0.
\end{align}
Let $s>\frac{5}{2}$, there exists a sufficiently small parameter $\delta$ such that, if $(w_{0},z_{0},\psi_{0})$ satisfies
\begin{align*}
&\Big[\|(w_0,z_0,\psi_0)\|_{H^{s}}
+C
\Big(\varepsilon\big(\|(m_{0},n_{0},r_{0})\|^{2}_{H^{s+2}}+\|(n_{0},r_0)\|_{H^{s+2}}\big)\Big)
\Big]\\
&\ \ \ \ \ \ \ \ \ \ \ \ \ \ \times e^{C\big(\|(m_{0},n_{0},r_{0})\|_{H^{s+2}}
+\|(n_{0},r_0)\|^{2}_{H^{s+2}}\big)}\leq\delta\min\{\nu,\eta,\lambda\},
\end{align*}
then the system (\ref{1.1}) has a unique global solution.
\end{thm}

\begin{rem}\label{rem1.1}
The initial data $(u_0, v_0, \theta_0)$ in Theorem \ref{thm1.1} is not small. In fact,
\begin{align*}
\|u_{0}\|_{H^{s}}&\geq\|U_{0}\|_{H^{s}}-\|w_{0}\|_{H^{s}}\\
&=\Big[\int_{\mr^3}(1+|\xi|^{2})^{s}|\xi|^{2}|\widehat{m}_{0}(\xi)|^{2}\md \xi\Big]^{\frac{1}{2}}-\|w_{0}\|_{H^{s}}\\
&\geq\Big[\int_{\mathcal{C}_{1}}(1+|\xi|^{2})^{s}|\xi|^{2}|\widehat{m }_{0}(\xi)|^{2}\md \xi\Big]^{\frac{1}{2}}-\|w_{0}\|_{H^{s}}\\
&\geq C\Big(\varepsilon^{-1}\log\log\frac{1}{\varepsilon}\Big)\varepsilon-\|w_{0}\|_{H^{s}}\\
&\geq C\log\log\frac{1}{\varepsilon}-\delta\min\{\nu,\eta,\lambda\},
\end{align*}
which will be really large when $\varepsilon$ is taken to be small. By applying the same argument, we can show that $\|v_{0}\|_{H^{s}}$ and $\|\theta_{0}\|_{H^{s}}$ are not small.
\end{rem}

When $d=2$, we can also obtain a unique global solution for the TCM. Likewise, we define three scalar functions $\overline{m}_0$, $\overline{n}_0$, $\overline{r}_0$ $\in C^{\infty}_{0}(\mr^2)$ with their Fourier transform satisfying
\begin{align}
&\widehat{\overline{m}}_0=
\Big(\varepsilon^{-\frac{1}{2}}\log\log\frac{1}{\varepsilon}\Big)\chi^*(\xi),&& \xi\in\mr^2 \label{1.6a},\\
&\widehat{\overline{n}}_0=\Big(\varepsilon^{-\frac{1}{2}}
\Big(\log\log\frac{1}{\varepsilon}\Big)^{\frac{1}{2}}\Big)\chi^*(\xi),&& \xi\in\mr^2 \label{1.7a},\\
&\widehat{\overline{r}}_0=\Big(\varepsilon^{-1}
\Big(\log\log\frac{1}{\varepsilon}\Big)^{\frac{1}{2}}\Big)\overline{\chi^*}(\xi),&& \xi\in\mr^2 \label{1.8a}
\end{align}
the smooth cutoff function $\chi^*(\xi)$ satisfies
\begin{center}
supp $\chi^*\subset\mathcal{D}$,$ \ $ $\chi^*\in[0,1]$ $ \ $ and $ \ $ $\chi^*=1$ $ \ $ on $ \ $ $\mathcal{D}_{1}$,\\
supp $\overline{\chi^*}\subset\mathcal{E}$,$ \ $ $\overline{\chi^*}\in[0,1]$ $ \ $ and $ \ $ $\overline{\chi^*}=1$ $ \ $ on $ \ $ $\mathcal{E}_{1}$,
\end{center}
where $\mathcal{D}$ and $\mathcal{D}_1$ denote the annuli,
\begin{align*}
&\mathcal{D}:=\Big\{\xi\in\mr^2|\ |\xi_{1}+\xi_{2}|\leq\varepsilon,1\leq|\xi|\leq 2\Big\},\\
&\mathcal{D}_{1}:=\Big\{\xi\in\mr^2|\ |\xi_{1}+\xi_{2}|\leq\varepsilon,\frac{4}{3}\leq|\xi|\leq\frac{5}{3}\Big\},\\
&\mathcal{E}:=\Big\{\xi\in\mr^2|\ |\xi_{1}+\xi_{2}|\leq\varepsilon,\varepsilon\leq|\xi|\leq 2\varepsilon\Big\},\\
&\mathcal{E}_{1}:=\Big\{\xi\in\mr^2|\ |\xi_{1}+\xi_{2}|\leq\varepsilon,\frac{4}{3}\varepsilon\leq|\xi|\leq\frac{5}{3}\varepsilon\Big\}.
\end{align*}

The result for 2-dimensional case can be stated as follows.
\begin{thm}\label{thm1.2}
Assume $\overline{m}_{0},\overline{n}_{0},\overline{r}_0$ are defined by (\ref{1.6a}), (\ref{1.7a}) and (\ref{1.8a}). Consider the initial data in the TCM system (\ref{1.1}) fulfills $\textup{div}w_{0}=0$ and
\begin{align*}
u_{0}=\overline{U}_{0}+w_{0},\ \ \ v_{0}=\overline{V}_{0}+z_{0}\ \ \ and \ \ \ \theta_0=\overline{\Theta}_0+\psi_0,
\end{align*}
where
\begin{align}\label{1.8}
\overline{U}_{0}=\nabla^\perp\overline{m}_0,\ \ \ \overline{V}_{0}=(\overline{n}_0, \overline{n}_0)\ \ \ and\ \ \ \overline{\Theta}_0=\overline{r}_0.
\end{align}
Let $s>2$, there exists a sufficiently small parameter $\delta$ such that, if $(w_{0},z_{0},\psi_{0})$ satisfies
\begin{align*}
&\Big[\|(w_0,z_0,\psi_0)\|_{H^{s}}
+C
\Big(\varepsilon\big(\|(\overline{m}_{0},\overline{n}_{0},\overline{r}_{0})\|^{2}_{H^{s+2}}
+\|(\overline{n}_{0},\overline{r}_{0})\|_{H^{s+2}}\big))
\Big]\\
&\ \ \ \ \ \ \ \ \ \ \ \ \ \ \times e^{C\big(\|(\overline{m}_{0},\overline{n}_{0},\overline{r}_{0})\|_{H^{s+2}}
+\|(\overline{n}_{0},\overline{r}_0)\|^{2}_{H^{s+2}}\big)}\leq\delta\min\{\nu,\eta,\lambda\},
\end{align*}
then the system (\ref{1.1}) has a unique global solution.
\end{thm}

\begin{rem}\label{rem1.2}
Li-Yu \cite{LY2020} studied the 2-dimensional TCM in (\ref{1.1}) with $\lambda=0$ and constructed a class of global solutions which permit the initial data $(u_0, v_0)$ large in $L^\infty(\mr^2)$, but the initial data $\theta_0$ is small in $L^\infty(\mr^2)$. In Theorem \ref{thm1.2}, we set the initial data $(u_0, v_0, \theta_0)$ can be large in $H^s(\mr^2)$ for $s>2$, which further improve the result in \cite{LY2020}.
\end{rem}

\section{Preliminaries}
\setcounter{equation}{0}
\vskip .1in

To prove Theorem \ref{thm1.1}, we renormalize the system (\ref{1.1}). Let $(m, n, r)$ be the solution of the following system
\begin{equation}\label{1.5}
\begin{cases}
{ \begin{array}{ll} \partial_{t}m+\nu m=0,\\
\partial_{t}n-\eta\Delta n+\partial_1 r+\partial_2 r+\partial_3 r=0,\\
\partial_{t}r+\lambda r=0,\\
m(x,0)=m_{0}(x),\ \ n(x,0)=n_{0}(x),\ \ r(x,0)=r_0(x),
 \end{array} }
\end{cases}
\end{equation}
which implies
\begin{align*}
m(t)=e^{-\nu t}m_0,\ \ \ n(t)=e^{\eta\Delta t}n_0-\int^t_0e^{\eta\Delta(t-\tau)}e^{-\lambda \tau}(\partial_1+\partial_2+\partial_3) r_0\md \tau,\ \ \ r(t)=e^{-\lambda t}r_0.
\end{align*}
Defining
\begin{center}
$U=\nabla\times m$, $\ \ $ $V=(n, n, n)$ $ \ \ $ and $ \ \ $ $\Theta=r$,
\end{center}
we can deduce from (\ref{1.5}) that
\begin{equation}\label{1.6}
\begin{cases}
{ \begin{array}{ll} \partial_{t}U+\nu U=0,\\
\partial_{t}V-\eta\Delta V+\nabla \Theta+\mathbb{A}\Theta=0,\\
\partial_{t}\Theta+\lambda\Theta=0,\\
U(x,0)=U_{0}(x),\ \ V(x,0)=V_{0}(x),\ \ \Theta(x,0)=\Theta_0(x),
 \end{array} }
\end{cases}
\end{equation}
where operator $\mathbb{A}=(\partial_2+\partial_3, \partial_1+\partial_3, \partial_1+\partial_2)$.\\
By (\ref{1.6}) and the definition of $U_0$, $V_0$, $\Theta_0$ in Theorem \ref{thm1.1}, $U$, $V$, $\Theta$ can be written as
\begin{align}\label{1.7}
\nonumber U(t)&=e^{-\nu t}U_{0}=e^{-\nu t}\nabla\times m_{0},\\
\nonumber V(t)&=e^{\eta\Delta t}V_{0}-\int^t_0 e^{\eta\Delta(t-\tau)}(\nabla \Theta+\mathbb{A}\Theta)\md \tau\\
\nonumber&=e^{\eta\Delta t}(n_{0}, n_{0}, n_{0})-\int^t_0 e^{\eta\Delta(t-\tau)}e^{-\lambda\tau}(\nabla r_0+\mathbb{A}r_0)\md\tau,\\
\Theta(t)&=e^{-\lambda t}\Theta_0=e^{-\lambda t}r_0.
\end{align}
Denoting $w=u-U$, $z=v-V$ and $\psi=\theta-\Theta$, then $(w,z,\psi)$ is the solution of the following equations
\begin{equation}\label{1.8}
\begin{cases}
{ \begin{array}{ll} \partial_{t}w+(w\cdot\nabla)w+(z\cdot\nabla)z+z(\nabla\cdot z)+\nu w+\nabla p\\
\ \ =f-(U\cdot\nabla)w-(w\cdot\nabla)U-(V\cdot\nabla)z-(z\cdot\nabla )V-z(\nabla\cdot V)-V(\nabla\cdot z),\\
\partial_{t}z+(w\cdot\nabla)z+(z\cdot\nabla)w
-\eta\Delta z+\nabla\psi\\
\ \ =g-(U\cdot\nabla)z-(w\cdot\nabla)V-(V\cdot\nabla)w-(z\cdot\nabla)U,\\
\partial_{t}\psi+(w\cdot\nabla)\psi+\lambda\psi+\nabla\cdot z=h-(w\cdot\nabla)\Theta-(U\cdot\nabla)\psi,\\
\textup{div}w=0,\\
w(x,0)=w_{0}(x),\ \ z(x,0)=z_{0}(x),\ \ \psi(x,0)=\psi_{0}(x),
 \end{array} }
\end{cases}
\end{equation}
where
\begin{align}\label{1.9}
\nonumber f=&-(U\cdot\nabla)U-(V\cdot\nabla)V-V(\nabla\cdot V),\\ \nonumber g=&-(U\cdot\nabla)V-(V\cdot\nabla)U+\mathbb{A}\Theta,\\
h=&-\nabla\cdot V-(U\cdot\nabla)\Theta.
\end{align}

Next, we present two lemmas which will be used in the proof of Theorem \ref{thm1.1}. Lemma \ref{lem2.1} offers upper bounds on $(U, V, \Theta)$ and $(f, g, h)$.
\begin{lem}\label{lem2.1}
Assume $m_{0}$, $n_{0}$ and $r_0$ are defined by (\ref{1.2}) and (\ref{1.3}), $U(t)$, $V(t)$ and $\Theta(t)$ by (\ref{1.7}), and $f$, $g$ and $h$ by (\ref{1.9}). Then for any $s>\frac{5}{2}$, the following estimates hold
\begin{equation}\label{2.6}
\begin{split}
&\|U\|_{L^1_tH^{s+1}}\leq C\|m_{0}\|_{H^{s+2}},\\
&\|V\|_{L^1_tH^{s+1}}\leq C(\|n_{0}\|_{H^{s+2}}+\|r_0\|_{H^{s+2}}),\\
&\|\Theta\|_{L^1_tH^{s+1}}\leq C\|r_{0}\|_{H^{s+2}},\\
&\|f\|_{L^1_tH^{s}}\leq C\varepsilon (\|m_{0}\|^{2}_{H^{s+2}}+\|n_{0}\|^{2}_{H^{s+2}}+\|r_0\|^2_{H^{s+2}}),\\
&\|g\|_{L^1_tH^{s}}\leq C\varepsilon (\|m_{0}\|^{2}_{H^{s+2}}+\|n_{0}\|^{2}_{H^{s+2}}+\|r_0\|^2_{H^{s+2}})++C\ep \|r_0\|_{H^{s}}\\
&\|h\|_{L^1_tH^{s}}\leq C\varepsilon (\|n_{0}\|_{H^{s}}+\|r_0\|_{H^{s+1}}+\|m_0\|^2_{H^{s+1}}+\|r_0\|^2_{H^{s+1}}).
\end{split}
\end{equation}
\begin{proof}
It is easy to prove the first three terms of (\ref{2.6}) by (\ref{1.7}) and the definition of $(m_0, n_0, r_0)$, we shall omit the detailed steps. Now, we prove the estimates of $f$, $g$ and $h$. Our main idea is to rewrite the terms in $f$, $g$ and $h$ so that each term contains $\partial_{i}+\partial_{j}$$(i, j=1, 2, 3, i\neq j)$.

For simplicity, here we take the first component of each term in $f$ as an example, the rest can be obtained in the same way.
\begin{align*}
&U=\nabla\times m=(\partial_{2}m_{3}-\partial_{3}m_{2}, \partial_{3}m_{1}-\partial_{1}m_{3}, \partial_{1}m_{2}-\partial_{2}m_{1}),\\
&V=(n, n, n).
\end{align*}
Direct calculations show that
\begin{equation*}
\begin{split}
U\cdot\nabla
U^{1}&=\partial_{2}m_{3}\partial_{1}\partial_{2}m_{3}-\partial_{1}m_{3}\partial_{2}\partial_{2}m_{3}
+\partial_{1}m_{3}\partial_{2}\partial_{3}m_{2}-\partial_{2}m_{3}\partial_{1}\partial_{3}m_{2}\\
&+\partial_{3}m_{1}\partial_{2}\partial_{2}m_{3}-\partial_{2}m_{1}\partial_{3}\partial_{2}m_{3}
+\partial_{2}m_{1}\partial_{3}\partial_{3}m_{2}-\partial_{3}m_{1}\partial_{2}\partial_{3}m_{2}\\
&+\partial_{1}m_{2}\partial_{3}\partial_{2}m_{3}-\partial_{3}m_{2}\partial_{1}\partial_{2}m_{3}
-\partial_{1}m_{2}\partial_{3}\partial_{3}m_{2}+\partial_{3}m_{2}\partial_{1}\partial_{3}m_{2}\\
&=(\partial_{2}+\partial_{1})m_{3}\partial_{1}\partial_{2}m_{3}
-\partial_{1}m_{3}(\partial_{1}+\partial_{2})\partial_{2}m_{3}\\
&+(\partial_{1}+\partial_{2})m_{3}\partial_{2}\partial_{3}m_{2}
-\partial_{2}m_{3}(\partial_{2}+\partial_{1})\partial_{3}m_{2}\\
&+(\partial_{3}+\partial_{2})m_{1}\partial_{2}\partial_{2}m_{3}
-\partial_{2}m_{1}(\partial_{2}+\partial_{3})\partial_{2}m_{3}\\
&+(\partial_{2}+\partial_{3})m_{1}\partial_{3}\partial_{3}m_{2}
-\partial_{3}m_{1}(\partial_{3}+\partial_{2})\partial_{3}m_{2}\\
&+(\partial_{1}+\partial_{3})m_{2}\partial_{3}\partial_{2}m_{3}
-\partial_{3}m_{2}(\partial_{3}+\partial_{1})\partial_{2}m_{3}\\
&-(\partial_{1}+\partial_{3})m_{2}\partial_{3}\partial_{3}m_{2}
+\partial_{3}m_{2}(\partial_{3}+\partial_{1})\partial_{3}m_{2}.
\end{split}
\end{equation*}
Taking the $H^{s}$ norm yields
\begin{equation*}
\begin{split}
\|U\cdot\nabla U^{1}\|_{H^{s}}\leq &\ C\big(\|(\partial_{2}+\partial_{1})m\|_{H^{s}}\|m\|_{H^{s+2}}+\|m\|_{H^{s+1}}\|(\partial_{1}+\partial_{2})m\|_{H^{s+1}}\\
&+\|(\partial_{1}+\partial_{2})m\|_{H^{s}}\|m\|_{H^{s+2}}+\|m\|_{H^{s+1}}\|(\partial_{2}+\partial_{1})m\|_{H^{s+1}}\\
&+\|(\partial_{3}+\partial_{2})m\|_{H^{s}}\|m\|_{H^{s+2}}+\|m\|_{H^{s+1}}\|(\partial_{2}+\partial_{3})m\|_{H^{s+1}}\\
&+\|(\partial_{2}+\partial_{3})m\|_{H^{s}}\|m\|_{H^{s+2}}+\|m\|_{H^{s+1}}\|(\partial_{3}+\partial_{2})m\|_{H^{s+1}}\\
&+\|(\partial_{1}+\partial_{3})m\|_{H^{s}}\|m\|_{H^{s+2}}+\|m\|_{H^{s+1}}\|(\partial_{3}+\partial_{1})m\|_{H^{s+1}}\\
&+\|(\partial_{1}+\partial_{3})m\|_{H^{s}}\|m\|_{H^{s+2}}+\|m\|_{H^{s+1}}\|(\partial_{3}+\partial_{1})m\|_{H^{s+1}}\big)\\
\leq&\
C\big(\|(\partial_{i}+\partial_{j})m\|_{H^{s}}\|m\|_{H^{s+2}}+\|m\|_{H^{s+1}}\|(\partial_{i}+\partial_{j})m\|_{H^{s+1}}\big),
\end{split}
\end{equation*}
where $i, j$ in the last line are summed over $i, j= 1, 2, 3$ and $i\neq j$. Similarly,
\begin{equation*}
V\cdot\nabla V^{1}
=n\partial_{1}n+n\partial_{2}n+n\partial_{3}n
=n\frac{\partial_{1}+\partial_{2}}{2}n+n\frac{\partial_{2}
+\partial_{3}}{2}n+n\frac{\partial_{3}+\partial_{1}}{2}n=V^{1} (\nabla\cdot V).
\end{equation*}
Taking the $H^{s}$ norm yields
\begin{equation*}
\begin{split}
&\|V\cdot\nabla V^{1}\|_{H^{s}}=\|V^{1} (\nabla\cdot V)\|_{H^{s}}\\
&\leq C\big(\|n\|_{H^{s}}\|(\partial_{1}+\partial_{2})n\|_{H^{s}}+\|n\|_{H^{s}}\|(\partial_{2}+\partial_{3})n\|_{H^{s}}
+\|n\|_{H^{s}}\|(\partial_{3}+\partial_{1})n\|_{H^{s}}\big)\\
&\leq C\varepsilon \big(e^{-2C_1t}\|n_0\|^2_{H^{s}}+\big(e^{-C_1(\cdot)}\ast e^{-C_0(\cdot)}\big)(t)\|r_0\|^2_{H^{s+1}}\big).
\end{split}
\end{equation*}
Where $\ast$ is convolution operator and $\big(e^{-C_1(\cdot)}\ast e^{-C_0(\cdot)}\big)(t)=\int^t_0 e^{-C_1(t-\tau)}e^{-C_0 \tau}\md \tau$.

Similar argument is used to the rest components, it holds that
\begin{equation*}
\begin{split}
\|f\|_{H^{s}}\leq &C\big(\|(\partial_{i}+\partial_{j})m\|_{H^{s}}\|m\|_{H^{s+2}}+\|m\|_{H^{s+1}}
\|(\partial_{i}+\partial_{j})m\|_{H^{s+1}}
\big)\\
&+C\varepsilon \big(e^{-2C_1t}\|n_0\|^2_{H^{s}}+\big(e^{-C_1(\cdot)}\ast e^{-C_0(\cdot)}\big)(t)\|r_0\|^2_{H^{s+1}}\big)\\
\leq &Ce^{-2C_{0}t}\big(\|(\partial_{i}+\partial_{j})m_{0}\|_{H^{s}}\|m_{0}\|_{H^{s+2}}
+\|m_{0}\|_{H^{s+1}}\|(\partial_{i}+\partial_{j})m_{0}\|_{H^{s+1}}
\big)\\
&+C\varepsilon \big(e^{-2C_1t}\|n_0\|^2_{H^{s}}+\big(e^{-C_1(\cdot)}\ast e^{-C_0(\cdot)}\big)(t)\|r_0\|^2_{H^{s+1}}\big)\\
\leq &C\varepsilon e^{-2C_{0}t}\|m_{0}\|^2_{H^{s+2}}+C\varepsilon \big(e^{-2C_1t}\|n_0\|^2_{H^{s}}+\big(e^{-C_1(\cdot)}\ast e^{-C_0(\cdot)}\big)(t)\|r_0\|^2_{H^{s+1}}\big)\\
\end{split}
\end{equation*}
Integrating both sides of the above inequality over $[0,t]$, by convolution Young's inequality, we deduce that
\begin{align*}
\|f\|_{L^1_tH^{s}}\leq& C\ep\|m_0\|^2_{H^{s+2}}\int^t_0e^{-2C_0\tau}\md \tau
+C\ep\|n_0\|^2_{H^{s+2}}\int^t_0e^{-2C_1\tau}\md \tau\\
&+C\ep\|r_0\|^2_{H^{s+2}}\int^t_0e^{-C_0\tau}\md \tau\int^t_0e^{-C_1\tau}\md \tau\\
\leq& C\ep\big(\|m_0\|^2_{H^{s+2}}+\|n_0\|^2_{H^{s+2}}+\|r_0\|^2_{H^{s+2}}\big).
\end{align*}
Since $\mathbb{A}=(\partial_2+\partial_3,\partial_1+\partial_3,\partial_1+\partial_3)$, we have
\begin{align*}
\|\mathbb{A}r\|_{H^s}\leq& \|(\partial_2+\partial_3)r\|_{H^s}+\|(\partial_1+\partial_3)r\|_{H^s}+
\|(\partial_1+\partial_3)r\|_{H^s}\\
\leq& C\ep e^{-C_0t}\|r_0\|_{H^s}.
\end{align*}
Then, by the similarly arguments that used in evaluating $\|f\|_{L^1_tH^s}$, $\|g\|_{L^1_tH^{s}}$ is estimated as
\begin{align}
\nonumber \|g\|_{L^1_tH^{s}}\leq C\varepsilon \big(\|m_0\|^2_{H^{s+2}}+\|n_0\|^2_{H^{s+2}}+\|r_0\|^2_{H^{s+2}}\big)+C\ep \|r_0\|_{H^{s}}.
\end{align}
For $h$, since
\begin{align*}
\|\nabla\cdot V\|_{H^s}=&\|(\partial_{1}+\partial_{2}+\partial_{3})n\|_{H^{s}}\\
=&\frac{1}{2}\|(\partial_{1}+\partial_{2})n+(\partial_{2}+\partial_{3})n+(\partial_{3}+\partial_{1})n\|_{H^{s}}\\
\leq& C\varepsilon e^{-C_{1}t}\|n_{0}\|_{H^{s}}+C\ep \big(e^{-C_1(\cdot)}\ast e^{-C_0(\cdot)}\big)(t)\|r_0\|_{H^{s+1}},
\end{align*}
and
\begin{align*}
\|(U\cdot\nabla)\Theta\|_{H^s}=&\|(\partial_2m_3-\partial_3m_2)\partial_1r
+(\partial_3m_1-\partial_1m_3)\partial_2r
+(\partial_1m_2-\partial_2m_1)\partial_3r\|_{H^s}\\
\leq&\|(\partial_2+\partial_1)m_3\partial_1r-\partial_1m_3(\partial_1+\partial_2)r\|_{H^s}\\
&+\|(\partial_3+\partial_2)m_1\partial_2r-\partial_2m_1(\partial_2+\partial_3)r\|_{H^s}\\
&+\|(\partial_1+\partial_3)m_2\partial_3r-\partial_3m_2(\partial_3+\partial_1)r\|_{H^s}\\
\leq& C\varepsilon(\|m\|_{H^s}\|r\|_{H^{s+1}}+\|m\|_{H^{s+1}}\|r\|_{H^s})\\
\leq& C\varepsilon e^{-2C_0t}(\|m_0\|^2_{H^{s+1}}+\|r_0\|^2_{H^{s+1}}),
\end{align*}
we have
\begin{align*}
\|h\|_{L^1_tH^{s}}\leq& C\varepsilon\|n_{0}\|_{H^{s}}\int^t_0 e^{-C_{1}\tau}\md\tau+C\ep\|r_0\|_{H^{s+1}}\int^t_0e^{-C_0\tau}\md \tau\int^t_0e^{-C_1\tau}\md \tau\\
&+C\varepsilon (\|m_0\|^2_{H^{s+1}}+\|r_0\|^2_{H^{s+1}})\int^t_0e^{-2C_0t}\md \tau\\
\leq& C\ep(\|n_0\|_{H^s}+\|r_0\|_{H^{s+1}}+\|m_0\|^2_{H^{s+1}}+\|r_0\|^2_{H^{s+1}}).
\end{align*}
This proves (\ref{2.6}).
\end{proof}
\end{lem}

Next, we recall the following commutator and bilinear estimates, the detailed proof can be found in \cite{KP1988}.
\begin{lem}\label{lem2.2} Let $s>0$. Let $p, p_{1}, p_{3}\in(1,\infty)$ and $p_{2}, p_{4}\in[1,\infty]$ satisfies
\begin{align*}
\frac{1}{p}=\frac{1}{p_{1}}+\frac{1}{p_{2}}=\frac{1}{p_{3}}+\frac{1}{p_{4}}.
\end{align*}
Then there exists constants $C$ such that
\begin{align*}
&\|[J^{s},F]G\|_{L^{p}}\leq\ C(\|J^{s}F\|_{L^{p_{1}}}\|G\|_{L^{p_{2}}}+\|J^{s-1}G\|_{L^{p_{3}}}\|\nabla F\|_{L^{p_{4}}}),\\
&\|J^{s}(FG)\|_{L^{p}}\leq\ C(\|J^{s}F\|_{L^{p_{1}}}\|G\|_{L^{p_{2}}}+\|J^{s}G\|_{L^{p_{3}}}\|F\|_{L^{p_{4}}}),
\end{align*}
where $J=(I-\Delta)^{\frac{1}{2}}$ and the commutator $[J^s,F]G=J^s(FG)-F(J^sG)$. The operator $J$ is defined via the Fourier transform
\begin{align*}
\widehat{Jf}(\xi)=(1+|\xi|^2)^{\frac{1}{2}}\widehat{f}(\xi),
\end{align*}
and
\begin{align*}
\widehat{f}(\xi)=\int_{\mr^d}e^{-ix\cdot\xi}f(x)\md x.
\end{align*}
\end{lem}

\section{Proof of Theorem \ref{thm1.1}}
\setcounter{equation}{0}
\vskip .1in

This section is devoted to the proof of Theorem \ref{thm1.1}. We derive the $H^{s}$ estimates of $(w, z, \psi)$ and apply the bootstrap argument to complete the proof.

\begin{proof}
Applying $J^{s}$ to $(\ref{1.8})_{1,2,3}$ and dotting with $(J^{s}w, J^{s}z, J^{s}\psi)$ yields
\begin{align}\label{3.1}
\frac{1}{2}\frac{\md}{\md t}\|(w, z, \psi)\|^{2}_{H^{s}}
+\nu\|w\|^2_{H^{s}}+\eta\|z\|^2_{H^{s+1}}+\lambda\|\psi\|^2_{H^{s}}
=\sum^{7}_{i=1}I_{i},
\end{align}
with
\begin{align*}
&I_{1}=-\int_{\mr^3}[J^{s},w\cdot\nabla]w\cdot J^{s}w\md x-\int_{\mr^3}[J^{s},w\cdot\nabla]z\cdot J^{s}z\md x-\int_{\mr^3}[J^{s},w\cdot\nabla]\psi\cdot J^{s}\psi\md x,\\
&I_{2}=-\int_{\mr^3}[J^{s},U\cdot\nabla]w\cdot J^{s}w\md x-\int_{\mr^3}[J^{s},U\cdot\nabla]z\cdot J^{s}z\md x-\int_{\mr^3}[J^{s},U\cdot\nabla]\psi\cdot J^{s}\psi\md x,\\
&I_{3}=-\int_{\mr^3}J^s(z\cdot\nabla z)\cdot J^{s}w\md x
-\int_{\mr^3}J^s(z(\nabla\cdot z))\cdot J^{s}w\md x
-\int_{\mr^3}J^{s}(z\cdot\nabla w)\cdot J^{s}z\md x,\\
&I_{4}=-\int_{\mr^3}J^{s}(w\cdot\nabla U)\cdot J^{s}w\md x
-\int_{\mr^3}J^{s}(z\cdot\nabla U)\cdot J^{s}z\md x
-\int_{\mr^3}J^{s}(w\cdot\nabla\Theta)\cdot J^{s}\psi\md x,\\
&I_{5}=-\int_{\mr^3}J^{s}(V\cdot\nabla z)\cdot J^{s}w\md x
-\int_{\mr^3}J^{s}(V(\nabla\cdot z))\cdot J^{s}w\md x
-\int_{\mr^3}J^{s}(V\cdot\nabla w)\cdot J^{s}z\md x,\\
&I_{6}=-\int_{\mr^3}J^{s}(z\cdot\nabla V)\cdot J^{s}w\md x
-\int_{\mr^3}J^{s}(z(\nabla\cdot V))\cdot J^{s}w\md x
-\int_{\mr^3}J^{s}(w\cdot\nabla V)\cdot J^{s}z\md x,\\
&I_{7}=\int_{\mr^3}J^{s}f\cdot J^{s}w\md x
+\int_{\mr^3}J^{s}g\cdot J^{s}z\md x
+\int_{\mr^3}J^{s}h\cdot J^{s}\psi\md x,
\end{align*}
where we have used the simple fact
\begin{align*}
\int_{\mr^3}J^{s}\nabla\psi\cdot J^{s}z\md x+\int_{\mr^3}J^{s}\nabla\cdot zJ^{s}\psi\md x=0.
\end{align*}
We now estimate the terms on the right hand side of (\ref{3.1}). By H\"{o}lder's inequality, Lemma \ref{lem2.2} and Young inequality, $I_{1}$ can be bounded by
\begin{equation*}
\begin{split}
|I_{1}|\leq &C\|[J^{s},w\cdot\nabla]w\|_{L^2}\|w\|_{H^{s}}
+C\|[J^{s},w\cdot\nabla]z\|_{L^2}\|z\|_{H^{s}}
+C\|[J^{s},w\cdot\nabla]\psi\|_{L^2}\|\psi\|_{H^{s}}\\
\leq& C\|\nabla w\|_{L^{\infty}}\|w\|^{2}_{H^{s}}
+C\|\nabla z\|_{L^{\infty}}\|w\|_{H^{s}}\|z\|_{H^{s}}
+C\|\nabla w\|_{L^{\infty}}\|z\|^{2}_{H^{s}}\\
&+C\|\nabla\psi\|_{L^{\infty}}\|w\|_{H^{s}}\|\psi\|_{H^{s}}+C\|\nabla w\|_{L^{\infty}}\|\psi\|^{2}_{H^{s}}\\
\leq& C\|w\|^{3}_{H^{s}}+C\|z\|_{H^{s+1}}\|w\|_{H^{s}}\|z\|_{H^{s}}
+C\|w\|_{H^{s}}\|\psi\|^{2}_{H^{s}}\\
\leq& C\|(w, z, \psi)\|_{H^{s}}(\|w\|^{2}_{H^{s}}+\|z\|^{2}_{H^{s+1}}
+\|\psi\|^{2}_{H^{s}}),\\
\end{split}
\end{equation*}
where we have used the following facts, for $s>\frac{5}{2}$
\begin{center}
$\|\nabla w\|_{L^{\infty}}\leq C\|w\|_{H^{s}}$,
$\|\nabla\psi\|_{L^{\infty}}\leq C\|\psi\|_{H^{s}}$ and
$\|\nabla z\|_{L^{\infty}}\leq C\|z\|_{H^{s}}\leq C\|z\|_{H^{s+1}}$.
\end{center}
For $I_{2}$, by the estimates of $U$ and $\nabla U$ in Lemma \ref{lem2.1}, we deduce
\begin{equation*}
\begin{split}
|I_{2}|\leq& C\|U\|_{H^{s}}(\|\nabla w\|_{L^{\infty}}\|w\|_{H^{s}}
+\|\nabla z\|_{L^{\infty}}\|z\|_{H^{s}}
+\|\nabla \psi\|_{L^{\infty}}\|\psi\|_{H^{s}})\\
&+C\|\nabla U\|_{L^{\infty}}
\|(w,z,\psi)\|^{2}_{H^{s}}\\
\leq& C\|U\|_{H^{s}}
\|(w,z,\psi)\|^{2}_{H^{s}}.
\end{split}
\end{equation*}
To evaluate $I_3$, we use the fact $z\cdot\nabla w=\nabla\cdot(w\otimes z)-w(\nabla\cdot z)$, then it holds that
\begin{align*}
|I_{3}|\leq& C\|z\cdot\nabla z\|_{H^{s}}\|w\|_{H^{s}}+C\|z(\nabla\cdot z)\|_{H^{s}}\|w\|_{H^s}+C\|w\otimes z\|_{H^{s}}\|\nabla z\|_{H^s}\\
&+C\|w(\nabla\cdot z)\|_{H^{s}}\|z\|_{H^s}\\
\leq& C\|z\|_{H^s}\|z\|_{H^{s+1}}\|w\|_{H^s}
+C\|w\|_{H^s}\|z\|_{H^s}\|z\|_{H^{s+1}}\\
\leq& C\|z\|_{H^s}(\|w\|^2_{H^s}+\|z\|^2_{H^{s+1}}).
\end{align*}
For $I_4$
\begin{align*}
|I_{4}|\leq& C\|w\cdot\nabla U\|_{H^s}\|w\|_{H^s}+C\|z\cdot\nabla U\|_{H^s}\|z\|_{H^s}+\|w\cdot\nabla\Theta\|_{H^s}\|\psi\|_{H^s}\\
\leq& C\|w\|_{H^{s}}\|\nabla U\|_{H^{s}}\|w\|_{H^{s}}
+C\|z\|_{H^{s}}\|\nabla U\|_{H^s}\|z\|_{H^{s}}
+C\|w\|_{H^{s}}\|\nabla\Theta\|_{H^{s}}\|\psi\|_{H^{s}}\\
\leq& C(\|U\|_{H^{s+1}}+\|\Theta\|_{H^{s+1}})
\|(w,z,\psi)\|^{2}_{H^{s}}.
\end{align*}
$I_{5}$ can be bounded by the similar idea of $I_3$
\begin{align*}
|I_{5}|\leq& C\|V\cdot\nabla z\|_{H^s}\|w\|_{H^{s}}
+C\|V(\nabla\cdot z)\|_{H^{s}}\|w\|_{H^{s}}
+C\|w\otimes V\|_{H^{s}}\|\nabla z\|_{H^s}\\
&+C\|w(\nabla\cdot V)\|_{H^{s}}\|z\|_{H^s}\\
\leq& C(\|V\|^2_{H^{s}}+\|V\|^2_{H^{s+1}})\|w\|^2_{H^s}+\frac{\eta}{2}\|z\|^2_{H^{s+1}}.
\end{align*}
By an argument similar to that used in evaluating $I_4$, it follows that
\begin{align*}
|I_{6}|\leq C\|V\|_{H^{s+1}}(\|w\|^{2}_{H^{s}}+\|z\|^{2}_{H^{s}}).
\end{align*}
Thanks to Lemma \ref{lem2.1}, the term $I_{7}$ can be bounded by
\begin{align*}
|I_{7}|\leq& \|f\|_{H^{s}}\|w\|_{H^{s}}+\|g\|_{H^{s}}\|z\|_{H^{s}}
+\|h\|_{H^{s}}\|\psi\|_{H^{s}}.\\
\end{align*}
Inserting all the estimates above for $I_{1}$ through $I_{7}$ in (\ref{3.1}) yields
\begin{align*}
\nonumber&\frac{1}{2}\frac{\md}{\md t}\|(w, z, \psi)\|^{2}_{H^{s}}
+\nu\|w\|^2_{H^{s}}+\frac{\eta}{2}\|z\|^2_{H^{s+1}}
+\lambda\|\psi\|^2_{H^{s}}\\
\nonumber\leq& C_{2}\|(w, z, \psi)\|_{H^{s}}(\|w\|^2_{H^{s}}+\|z\|^2_{H^{s+1}}+\|\psi\|^2_{H^{s}})\\ &+C_{3}(\|U\|_{H^{s+1}}+\|V\|_{H^{s+1}}+\|\Theta\|_{H^{s+1}}+\|V\|^2_{H^{s+1}})
\|(w, z, \psi)\|^{2}_{H^{s}}\\
&+C_{4}(\|f\|_{H^s}+\|g\|_{H^s}+\|h\|_{H^s})
\|(w, z, \psi)\|_{H^{s}}.
\end{align*}
Then we deduce
\begin{align}\label{3.5}
\nonumber&\frac{\md}{\md t}\|(w, z, \psi)\|^{2}_{H^{s}}
+\big(\nu-2C_{2}\|(w, z, \psi)\|_{H^{s}})\|w\|^2_{H^{s}}\\
\nonumber&+\big(\eta-2C_{2}\|(w, z, \psi)\|_{H^{s}}\big)\|z\|^2_{H^{s+1}}
+\big(\lambda-2C_{2}\|(w, z, \psi)\|_{H^{s}})\|\psi\|^2_{H^{s}}\\
\nonumber\leq& 2C_{3}(\|U\|_{H^{s+1}}+\|V\|_{H^{s+1}}+\|\Theta\|_{H^{s+1}}+\|V\|^2_{H^{s+1}})
\|(w, z, \psi)\|^{2}_{H^{s}}\\
&+2C_{4}(\|f\|_{H^s}+\|g\|_{H^s}+\|h\|_{H^s})
\|(w, z, \psi)\|_{H^{s}}.
\end{align}

Assume that $\|(w, z, \psi)\|_{H^{s}}$ is bounded, we prove $\|(w, z, \psi)\|_{H^{s}}$ actually admits a smaller bound when $\|(w_{0}, z_{0}, \psi_{0})\|_{H^{s}}$ is taken to be sufficiently small. To apply the bootstrap argument to (\ref{3.5}), we make the ansatz
\begin{align}\label{3.6}
\|(w, z, \psi)\|_{H^{s}}\leq\ M:=\frac{1}{2C_{2}}
\min\{\nu, \eta, \lambda\}.
\end{align}
Obviously, (\ref{3.6}) implies
\begin{align*}
\nu-2C_{2}\|(w, z, \psi)\|_{H^{s}}\geq 0
,\ \ \ \ \eta-2C_{2}\|(w, z, \psi)\|_{H^{s}}\geq 0,\ \ \ \ \lambda-2C_{2}\|(w, z, \psi)\|_{H^{s}}\geq 0.
\end{align*}
Thus we observe that (\ref{3.5}) leads to
\begin{equation*}
\begin{split}
\frac{\md}{\md t}&\|(w, z, \psi)\|_{H^{s}}\\
\leq& 2C_{3}(\|U\|_{H^{s+1}}+\|V\|_{H^{s+1}}+\|\Theta\|_{H^{s+1}}+\|V\|^2_{H^{s+1}})
\|(w, z, \psi)\|_{H^{s}}\\
&+2C_{4}(\|f\|_{H^s}+\|g\|_{H^s}+\|h\|_{H^s}).
\end{split}
\end{equation*}
By Gr\"{o}nwall inequality, we obtain
\begin{align}\label{3.7}
\nonumber&\|(w, z, \psi)\|_{H^{s}}\\
\nonumber\leq& e^{2C_{3}
\int^{t}_{0}(\|U\|_{H^{s+1}}+\|V\|_{H^{s+1}}+\|\Theta\|_{H^{s+1}}+\|V\|^2_{H^{s+1}})\md\tau}
\Big[\|(w_0, z_0, \psi_0)\|_{H^{s}}\\
&\nonumber+2C_{4}
\int^{t}_{0}(\|f\|_{H^s}+\|g\|_{H^s}+\|h\|_{H^s})\md\tau\Big]\\
\nonumber\leq& e^{C\big(\|(m_{0},n_{0},r_{0})\|_{H^{s+2}}
+\|(n_{0},r_0)\|^{2}_{H^{s+2}}\big)}
\Big[\|(w_0, z_0, \psi_0)\|_{H^{s}}\\
\nonumber&+C
\varepsilon\big(\|(m_{0},n_0,r_0)\|^{2}_{H^{s+2}}+\|(n_{0},r_0)\|_{H^{s+2}}\big)
\Big]\\
\leq& \delta\min\{\nu,\eta,\lambda\}.
\end{align}
Here we take
\begin{center}
$\delta=\frac{1}{4C_{2}},$
\end{center}
$\|(w, z, \psi)\|_{H^{s}}$ admits a smaller upper bound
\begin{center}
$\|(w, z, \psi)\|_{H^{s}}\leq\frac{M}{2},$
\end{center}
the bootstrap argument then implies that
\begin{center}
$\|(w, z, \psi)\|_{H^{s}}
\leq\frac{1}{4C_3}\min\{\nu,\eta,\lambda\}$\ \ \ \ \ for  $0<t<\infty$.
\end{center}
This completes the proof of Theorem \ref{thm1.1}.
\end{proof}

\section{Proof of Theorem \ref{thm1.2}}
In this section, we shall just provide the preparations of the proof of Theorem \ref{thm1.2}, the main proof follows the same argument of that for Theorem \ref{thm1.1}.

As a preparation of the proof, we first renormalize the system (\ref{1.1}) in the 2D case and provide several global upper bounds.
To this end, we assume that $(\overline{m},\overline{n},\overline{r})$ is the solution of the following system
\be\label{4.1}
\begin{cases}
{\begin{array}{ll}
\partial_t\overline{m}+\nu\overline{m}=0,\\
\partial_t\overline{n}-\eta\Delta\overline{n}+\partial_1\overline{r}+\partial_2\overline{r}=0,\\
\partial_t\overline{r}+\lambda\overline{r}=0.\\
\end{array}}
\end{cases}
\ee
Defining
\begin{align}\label{4.2}
\overline{U}=\nabla^\perp \overline{m},\ \ \ \ \ \ \overline{V}=(\overline{n},\overline{n}),\ \ \ \ \ \ \overline{\Theta}=\overline{r}.
\end{align}
Taking the differences $w=u-\overline{U},z=v-\overline{V},\psi=\theta-\overline{\Theta}$, it is easy to see that $(w,z,\psi)$ satisfy (\ref{1.8}) with $U, V,\Theta$ and $f,g,h$ replaced by $\overline{U},\overline{V},\overline{\Theta}$ and $\overline{f},\overline{g},\overline{h}$, respectively. Where $\overline{f},\overline{g},\overline{h}$ are defined as follows
\begin{align}\label{4.3}
\nonumber\overline{f}=&-(\overline{U}\cdot\nabla)\overline{U}-(\overline{V}\cdot\nabla)\overline{V}
-\overline{V}(\nabla\cdot\overline{V}),\\
\nonumber\overline{g}=&-(\overline{U}\cdot\nabla)\overline{V}-(\overline{V}\cdot\nabla)\overline{U}+\overline{\mathbb{A}}\ \overline{\Theta},\\
\overline{h}=&-\nabla\cdot\overline{V}-(\overline{U}\cdot\nabla)\overline{\Theta},
\end{align}
with $\overline{\mathbb{A}}=(\partial_2,\partial_1).$

Now we are in a position to prove upper bounds for $\overline{U},\overline{V},\overline{\Theta}$ and $\overline{f},\overline{g},\overline{h}$.
\begin{lem}\label{lem4.1}
Assume $\overline{m}_{0}$, $\overline{n}_{0}$ and $\overline{r}_0$ are defined by (\ref{1.6a}), (\ref{1.7a}) and (\ref{1.8a}), $\overline{U}(t)$, $\overline{V}(t)$ and $\overline{\Theta}(t)$ by (\ref{4.2}), and $\overline{f}$, $\overline{g}$ and $\overline{h}$ by (\ref{4.3}). Then for any $s>2$, the following estimates hold
\begin{equation}\label{4.4}
\begin{split}
&\|\overline{U}\|_{L^1_tH^{s+1}}\leq C\|\overline{m}_{0}\|_{H^{s+2}},\\
&\|\overline{V}\|_{L^1_tH^{s+1}}\leq C(\|\overline{n}_{0}\|_{H^{s+2}}+\|\overline{r}_{0}\|_{H^{s+2}}),\\
&\|\overline{\Theta}\|_{L^1_tH^{s+1}}\leq C\|\overline{r}_{0}\|_{H^{s+1}},\\
&\|\overline{f}\|_{L^1_tH^{s}}, \|\overline{g}\|_{L^1_tH^{s}}\leq C\varepsilon (\|\overline{m}_{0}\|^{2}_{H^{s+2}}+\|\overline{n}_{0}\|^{2}_{H^{s+2}}+\|\overline{r}_{0}\|^2_{H^{s+2}}),\\
&\|\overline{h}\|_{L^1_tH^{s}}\leq C\varepsilon (\|\overline{n}_{0}\|_{H^{s}}+\|\overline{r}_{0}\|_{H^{s+1}}+\|\overline{m}_{0}\|^2_{H^{s+1}}+\|\overline{r}_{0}\|^2_{H^{s+1}}).
\end{split}
\end{equation}
\begin{proof}
It is easy to prove the first three terms of (\ref{4.4}) by (\ref{4.2}) and the definition of $(\overline{m}_{0}, \overline{n}_{0}, \overline{r}_{0})$. To prove the upper bound for $\|f\|_{H^s}$, we rewrite the first component of $f$ as
\begin{align*}
\overline{f}_1=&-\partial_2\overline{m}\partial_1\partial_2\overline{m}
+\partial_1\overline{m}\partial_2\partial_2\overline{m}
-2\partial_1\overline{n}\partial_1\partial_1\overline{n}
+2\partial_2\overline{n}\partial_2\partial_1\overline{n}\\
=&-(\partial_1+\partial_2)\overline{m}\partial_1\partial_2\overline{m}
+\partial_1\overline{m}(\partial_1+\partial_2)\partial_2\overline{m}
-2(\partial_1+\partial_2)\overline{n}\partial_1\partial_1\overline{n}
+2\partial_2\overline{n}(\partial_1+\partial_2)\partial_1\overline{n}.
\end{align*}
Taking the $H^s$ norm, by H\"{o}lder inequality and Sobolev embedding, we deduce that
\begin{align*}
\|\overline{f}_1\|_{H^s}\leq& C\big(\|(\partial_1+\partial_2)\overline{m}\|_{H^s}\|\overline{m}\|_{H^{s+2}}
+\|(\partial_1+\partial_2)\overline{m}\|_{H^{s+1}}\|\overline{m}\|_{H^{s+1}}\\
&+\|(\partial_1+\partial_2)\overline{n}\|_{H^s}\|\overline{n}\|_{H^{s+2}}
+\|(\partial_1+\partial_2)\overline{n}\|_{H^{s+1}}\|\overline{n}\|_{H^{s+1}}\big)\\
\leq& C\ep e^{-2C_0t}\|\overline{m}_0\|^2_{H^{s+2}}+C\ep e^{-2C_1t}\|\overline{n}_0\|^2_{H^{s+2}}+C\ep \Big(e^{-C_1(\cdot)}\ast e^{-C_0(\cdot)}\Big)(t)\|\overline{r}_0\|^2_{H^{s+2}}.
\end{align*}
In addition, $\overline{f}_2$ admits the same bound, which implies that
\begin{align*}
\|\overline{f}\|_{L^1_tH^s}\leq C\varepsilon (\|\overline{m}_{0}\|^{2}_{H^{s+2}}+\|\overline{n}_{0}\|^{2}_{H^{s+2}}
+\|\overline{r}_{0}\|^2_{H^{s+2}}).
\end{align*}
The upper bounds for $\|\overline{f}\|_{L^1_tH^s},\|\overline{h}\|_{L^1_tH^s}$ can be similarly obtained. The proof of Lemma \ref{lem4.1} is thus complete.
\end{proof}
\end{lem}
The proof of Theorem \ref{thm1.2} is the same as that of Theorem \ref{thm1.1}, so we omit the details. Therefore, the proofs of our main result have been complete.

\textbf{Acknowledgements} The research of B Yuan
was partially supported by the National Natural Science Foundation
of China (No. 11471103).


\end{document}